\documentclass{amsart}
\usepackage{amsmath}
\usepackage{amsfonts}

\newtheorem{theorem}{Theorem}
\theoremstyle{plain}

\newtheorem{corollary}{Corollary}

\numberwithin{equation}{section}
\input{tcilatex}

\begin{document}
\title[$q$\textbf{-Genocchi polynomials}]{\textbf{Identities involving }$q$%
\textbf{-Genocchi numbers and polynomials}}
\author[\textbf{S. Arac\i }]{\textbf{Serkan Arac\i }}
\address{\textbf{University of Gaziantep, Faculty of Science and Arts,
Department of Mathematics, 27310 Gaziantep, TURKEY}}
\email{\textbf{mtsrkn@hotmail.com; saraci88@yahoo.com.tr}}
\author[\textbf{M. Acikgoz}]{\textbf{Mehmet Acikgoz}}
\address{\textbf{University of Gaziantep, Faculty of Science and Arts,
Department of Mathematics, 27310 Gaziantep, TURKEY}}
\email{\textbf{acikgoz@gantep.edu.tr}}
\author[\textbf{H. Jolany}]{\textbf{Hassan Jolany}}
\address{\textbf{School of Mathematics, Statistics and Computer Science,
University of Tehran, Iran}}
\email{\textbf{hassan.jolany@khayam.ut.ac.ir; hassan.jolany@mail.com}}
\author[\textbf{Y. He}]{\textbf{Yuan He}}
\address{\textbf{Department of Mathematics, Kunming University of Science
and Technology, Kunming, Yuannan 650500, People' Republic of China}}
\email{\textbf{hyyhe@yahoo.com.cn}}

\begin{abstract}
In this paper, we focus on the $q$-Genocchi numbers and polynomials. We
shall introduce new identities of the q-Genocchi numbers and polynomials by
using the fermionic $p$-adic integral on $%
\mathbb{Z}
_{p}$ which are very important in the study of Frobenius-Genocchi numbers
and polynomials. Also, we give Cauchy-integral formula for the $q$-Genocchi
polynomials and moreover by using measure theory on $p$-adic integral we
derive the distribution formula $q$-Genocchi polynomials. Finally, we
present a new definition of $q$-Zeta-type function by using Mellin
transformation which is the interpolation function of the $q$-Genocchi
polynomials at negative integers. 
\end{abstract}

\maketitle

\section{\textbf{PRELIMINARIES}}

Recently, \textit{Kim} and \textit{Lee} have given some properties for the $%
q $-Euler numbers and polynomials in \cite{Kim}. Actually, this type numbers
and polynomials and their $q$-extensions or variously generalizations have
been studied in several different ways for a long time (for details, see
[1-21]).

By using $p$-adic $q$-integral on $%
\mathbb{Z}
_{p}$, \textit{Kim} defined many new generating functions of the $q$%
-Bernoulli poynomials, $q$-Euler polynomials and $q$-Genocchi polynomials in
his arithmetic works (for details, see [6-19]). The works of \textit{Kim}
have been benefited for further works of many mathematicians in Analytic
numbers theory. Actually, we motivated from his inspiring works to write
this paper.

Assume that $p$ be a fixed odd prime number. Throughout this work, we
require the definitions of the some notations such that let $%
\mathbb{Q}
_{p}$ be the field $p$-adic rational numbers and let $%
\mathbb{C}
_{p}$ be the completion of algebraic closure of $%
\mathbb{Q}
_{p}$. That is, 
\begin{equation*}
\boldsymbol{%
\mathbb{Q}
}_{p}=\left\{ x=\sum_{n=-k}^{\infty }a_{n}p^{n}:0\leq a_{n}\leq p-1\right\} 
\text{.}
\end{equation*}

Then $%
\mathbb{Z}
_{p}$ is integral domain which is defined by 
\begin{equation*}
\boldsymbol{%
\mathbb{Z}
}_{p}=\left\{ x=\sum_{n=0}^{\infty }a_{n}p^{n}:0\leq a_{n}\leq p-1\right\}
\end{equation*}

or%
\begin{equation*}
\boldsymbol{%
\mathbb{Z}
}_{p}=\left\{ x\in 
\mathbb{Q}
_{p}:\left\vert x\right\vert _{p}\leq 1\right\} \text{.}
\end{equation*}

We assume that $q\in 
\mathbb{C}
_{p}$ with $\left\vert 1-q\right\vert _{p}<1$ as an indeterminate. The $p$%
-adic absolute value $\left\vert .\right\vert _{p}$, is normally given by 
\begin{equation*}
\left\vert x\right\vert _{p}=p^{-r}
\end{equation*}

where $x=p^{r}\frac{s}{t}$ with $\left( p,s\right) =\left( p,t\right)
=\left( s,t\right) =1$ and $r\in 
\mathbb{Q}
$.

The $q$-extension of $x$ with the display notation of $\left[ x\right] _{q}$
is introduced by 
\begin{equation*}
\left[ x\right] _{q}=\frac{1-q^{x}}{1-q}\text{.}
\end{equation*}%
We want to note that $\lim_{q\rightarrow 1}\left[ x\right] _{q}=x$\
(see[1-22]).

Also, we use notation $%
\mathbb{N}
^{\ast }$ for the union of zero and natural numbers.

We consider that $\eta $ is a uniformly differentiable function at a point $%
a\in 
\mathbb{Z}
_{p}$, if the difference quotient 
\begin{equation*}
\Phi _{\eta }\left( x,y\right) =\frac{\eta \left( x\right) -\eta \left(
y\right) }{x-y},
\end{equation*}%
have a limit $\eta 
{\acute{}}%
\left( a\right) $ as $\left( x,y\right) \rightarrow \left( a,a\right) $ and
denote this by $\eta \in UD\left( 
\mathbb{Z}
_{p}\right) $. \bigskip Then, for $\eta \in UD\left( 
\mathbb{Z}
_{p}\right) $, we can discuss the following%
\begin{equation*}
\frac{1}{\left[ p^{n}\right] _{q}}\sum_{0\leq \xi <p^{n}}\eta \left( \xi
\right) q^{\xi }=\sum_{0\leq \xi <p^{n}}\eta \left( \xi \right) \mu
_{q}\left( \xi +p^{n}%
\mathbb{Z}
_{p}\right) \text{,}
\end{equation*}

which represents as a $p$-adic $q$-analogue of Riemann sums for $\eta $. \
The integral of $\eta $ on $%
\mathbb{Z}
_{p}$ will be defined as the limit $\left( n\rightarrow \infty \right) $ of
these sums, when it exists. The $p$-adic $q$-integral of function $\eta \in
UD\left( 
\mathbb{Z}
_{p}\right) $ is defined by T. Kim in \cite{KIM}, \cite{Kim 3}, \cite{Kim 8}
by 
\begin{equation}
I_{q}\left( \eta \right) =\int_{%
\mathbb{Z}
_{p}}\eta \left( \xi \right) d\mu _{q}\left( \xi \right) =\lim_{n\rightarrow
\infty }\frac{1}{\left[ p^{n}\right] _{q}}\sum_{\xi =0}^{p^{n}-1}\eta \left(
\xi \right) q^{\xi }\text{.}  \label{equation 1}
\end{equation}

The bosonic integral is considered as a bosonic limit $q\rightarrow 1,$ $%
I_{1}\left( \eta \right) =\lim_{q\rightarrow 1}I_{q}\left( \eta \right) $.
Similarly, the fermionic $p$-adic integral on $%
\mathbb{Z}
_{p}$ is introduced by T. Kim as follows:%
\begin{equation}
I_{-q}\left( \eta \right) =\int_{%
\mathbb{Z}
_{p}}\eta \left( \xi \right) d\mu _{-q}\left( \xi \right)  \label{equation 2}
\end{equation}

(for more details, see [16-18]). Obviously that%
\begin{equation}
\lim_{q\rightarrow 1}I_{-q}\left( \eta \right) =I_{-1}\left( \eta \right)
=\int_{%
\mathbb{Z}
_{p}}\eta \left( \xi \right) d\mu _{-1}\left( \xi \right)
=\lim_{n\rightarrow \infty }\sum_{\xi =0}^{p^{n}-1}\eta \left( \xi \right)
\left( -1\right) ^{\xi }\text{.}  \label{equation 3}
\end{equation}

From (\ref{equation 3}), it is well-known as the useful property for the
fermionic $p$-adic $q$-integral on $%
\mathbb{Z}
_{p}$:%
\begin{equation}
I_{-1}\left( \eta _{1}\right) +I_{-1}\left( \eta \right) =2\eta \left(
0\right) \text{,}  \label{equation 10}
\end{equation}

where $\eta _{1}\left( x\right) =\eta \left( x+1\right) $ (for details,
see[2-4, 11-14, 16-22]).

The $q$-Genocchi polynomials can be introduced as follows:%
\begin{equation}
G_{n,q}\left( x\right) =n\int_{%
\mathbb{Z}
_{p}}q^{\xi }\left( x+\xi \right) ^{n-1}d\mu _{-1}\left( \xi \right) \text{.}
\label{equation 108}
\end{equation}

From (\ref{equation 108}), we have%
\begin{equation*}
G_{n,q}\left( x\right) =\sum_{l=0}^{n}\binom{n}{l}x^{l}G_{n-l,q}
\end{equation*}

where $G_{n,q}(0):=G_{n,q}$ are called $q$-Genocchi numbers. Then, $q$%
-Genocchi numbers can be given by%
\begin{equation*}
G_{0,q}=0\text{ and }q\left( G_{q}+1\right) ^{n}+G_{n,q}=\left\{ 
\begin{array}{cc}
2, & \text{if }n=1 \\ 
0, & \text{if\ }n\geq 1%
\end{array}%
\right.
\end{equation*}

with the usual convention about replacing $\left( G_{q}\right) ^{n}$ by $%
G_{n,q}$ is used (for details, see \cite{Araci 4}).

Our objective in the present paper is to derive not only new but also novel
and interesting properties of the $q$-Genocchi numbers and polynomials. Our
applications for the $q$-Genocchi polynomials seem to be useful in
mathematics for engineerings (on this subject, see \cite{Andrews}).

\section{\textbf{SOME PROPERTIES ON THE }$q$\textbf{-GENOCCHI NUMBERS AND
POLYNOMIALS}}

Let $\eta \left( x\right) =q^{x}e^{t\left( x+\xi \right) }$. Then, by using (%
\ref{equation 10}), we see that%
\begin{equation*}
t\int_{%
\mathbb{Z}
_{p}}q^{\xi }e^{t\left( x+\xi \right) }d\mu _{-1}\left( \xi \right) =\frac{2t%
}{qe^{t}+1}e^{xt}\text{.}
\end{equation*}

From the last equality and (\ref{equation 108}), we easily derive the
following generating function of the $q$-Genocchi polynomials: 
\begin{equation}
t\int_{%
\mathbb{Z}
_{p}}q^{\xi }e^{t\left( x+\xi \right) }d\mu _{-1}\left( \xi \right)
=\sum_{n=0}^{\infty }G_{n,q}\left( x\right) \frac{t^{n}}{n!}=\frac{2t}{%
qe^{t}+1}e^{xt}\text{, }\left\vert \log q+t\right\vert <\pi \text{.}
\label{equation 100}
\end{equation}

Substituting $x\rightarrow x+y$ into (\ref{equation 100}), then we write%
\begin{gather*}
\sum_{n=0}^{\infty }G_{n,q}\left( x+y\right) \frac{t^{n}}{n!}=\frac{2t}{%
qe^{t}+1}e^{\left( x+y\right) t} \\
=\left( \sum_{n=0}^{\infty }G_{n,q}\left( x\right) \frac{t^{n}}{n!}\right)
\left( \sum_{n=0}^{\infty }y^{n}\frac{t^{n}}{n!}\right) \\
=\sum_{n=0}^{\infty }\left( \sum_{k=0}^{n}\binom{n}{k}G_{k,q}\left( x\right)
y^{n-k}\right) \frac{t^{n}}{n!}
\end{gather*}

From the above, we easily express the following theorem:

\begin{theorem}
The following holds:%
\begin{equation}
G_{n,q}\left( x+y\right) =\sum_{k=0}^{n}\binom{n}{k}G_{k,q}\left( x\right)
y^{n-k}\text{.}  \label{equation 101}
\end{equation}
\end{theorem}

By (\ref{equation 101}), we consider the following%
\begin{equation*}
G_{n,q}\left( x+y\right) =\frac{2}{\left[ 2\right] _{q}}ny^{n-1}+%
\sum_{k=2}^{n}\binom{n}{k}G_{k,q}\left( x\right) y^{n-k}\text{.}
\end{equation*}

From this, it follows that%
\begin{equation*}
G_{n,q}\left( x+y\right) -\frac{2}{\left[ 2\right] _{q}}ny^{n-1}=%
\sum_{k=2}^{n}\binom{n}{k}G_{k,q}\left( x\right) y^{n-k}
\end{equation*}

can be derived and so we reach the following theorem:

\begin{theorem}
For $n\in 
\mathbb{N}
^{\ast }$, one has%
\begin{gather}
\sum_{k=0}^{n}\frac{\binom{n}{k}}{\left( k+2\right) \left( k+1\right) }%
G_{k+2,q}\left( x\right) y^{n-k}  \label{equation 102} \\
=\frac{G_{n+2,q}\left( x+y\right) -\frac{2}{\left[ 2\right] _{q}}\left(
n+2\right) y^{n+1}}{\left( n+2\right) \left( n+1\right) }\text{.}  \notag
\end{gather}
\end{theorem}

Replacing $y$ by $-y$ into (\ref{equation 102}), then we get%
\begin{gather}
\frac{G_{n+2,q}\left( x-y\right) -\frac{2}{\left[ 2\right] _{q}}\left(
n+2\right) \left( -1\right) ^{n+1}y^{n+1}}{\left( n+2\right) \left(
n+1\right) }  \label{equation 103} \\
=\sum_{k=0}^{n}\frac{\binom{n}{k}\left( -1\right) ^{n-k}}{\left( k+2\right)
\left( k+1\right) }G_{k+2,q}\left( x\right) y^{n-k}\text{.}  \notag
\end{gather}

By (\ref{equation 103}), it follows that%
\begin{gather}
\sum_{k=0}^{n}\frac{\binom{n}{k}\left( -1\right) ^{k}}{\left( k+2\right)
\left( k+1\right) }G_{k+2,q}\left( x\right) y^{n-k}  \label{equation 104} \\
=\frac{\left( -1\right) ^{n}G_{n+2,q}\left( x-y\right) +\frac{2}{\left[ 2%
\right] _{q}}\left( n+2\right) y^{n+1}}{\left( n+1\right) \left( n+2\right) }%
\text{.}  \notag
\end{gather}

Therefore, from the expressions of (\ref{equation 102}) and (\ref{equation
104}), we procure the following theorem:

\begin{theorem}
The following holds true:%
\begin{gather}
\sum_{k=0}^{\left[ \frac{n}{2}\right] }\frac{\binom{n}{2k}}{\left(
k+1\right) \left( 2k+1\right) }G_{2k+2,q}\left( x\right) y^{n-2k}
\label{equation 105} \\
=\frac{\left( -1\right) ^{n}G_{n+2,q}\left( x-y\right) +G_{n+2,q}\left(
x+y\right) }{\left( n+1\right) \left( n+2\right) }  \notag
\end{gather}%
where $\left[ .\right] $ is Gauss' symbol.
\end{theorem}

By (\ref{equation 103}), we have the following identity%
\begin{gather}
\sum_{k=2}^{n}\frac{\binom{n}{k}\left( -1\right) ^{k}}{k\left( k-1\right) }%
G_{k,q}\left( x\right) y^{n-k}  \label{equation 106} \\
=\frac{\left( -1\right) ^{n}G_{n,q}\left( x-y\right) +\frac{2}{\left[ 2%
\right] _{q}}ny^{n-1}}{n\left( n-1\right) }\text{.}  \notag
\end{gather}

By (\ref{equation 103}), (\ref{equation 104}) and (\ref{equation 106}), then
we have the following theorem:

\begin{theorem}
For $n\in 
\mathbb{N}
^{\ast }$, we get%
\begin{gather}
\frac{\left( -1\right) ^{n}G_{n+2,q}\left( x-y\right) +G_{n+2,q}\left(
x+y\right) }{\left( n+2\right) \left( n+1\right) }  \label{equation 107} \\
=\sum_{k=0}^{\left[ \frac{n+1}{2}\right] }\frac{\binom{n}{2k}}{\left(
k+1\right) \left( 2k+1\right) }G_{2k+2,q}\left( x\right) y^{n-2k}\text{.} 
\notag
\end{gather}
\end{theorem}

Taking $y=1$ into (\ref{equation 102}), it follows that%
\begin{equation}
q\sum_{k=0}^{n}\frac{\binom{n}{k}}{\left( k+2\right) \left( k+1\right) }%
G_{k+2,q}\left( x\right) =\frac{qG_{n+2,q}\left( x+1\right) }{\left(
n+1\right) \left( n+2\right) }-\frac{2q}{\left[ 2\right] _{q}\left(
n+1\right) }\text{.}  \label{equation 109}
\end{equation}

We need the following for sequel of this paper:%
\begin{equation*}
2e^{tx}=\frac{1}{t}\left( q\frac{2t}{qe^{t}+1}e^{\left( x+1\right) t}+\frac{%
2t}{qe^{t}+1}e^{xt}\right) \text{.}
\end{equation*}

From the above, we easily develop the following:%
\begin{equation}
qG_{n+1,q}\left( x+1\right) +G_{n+1,q}\left( x\right) =\left( n+1\right)
2x^{n}\text{.}  \label{equation 110}
\end{equation}

By (\ref{equation 109}) and (\ref{equation 110}), we state the following
theorem:

\begin{theorem}
The following holds:%
\begin{gather}
\sum_{k=0}^{n}\frac{\binom{n}{k}}{\left( k+2\right) \left( k+1\right) }%
G_{k+2,q}\left( x\right)  \label{equation 111} \\
=\frac{2x^{n+1}}{qn+q}-\frac{G_{n+2,q}\left( x\right) }{\left( qn+q\right)
\left( n+2\right) }-\frac{2}{\left[ 2\right] _{q}\left( n+1\right) }\text{.}
\notag
\end{gather}
\end{theorem}

Thanks to equality of $\lim_{q\rightarrow 1}G_{n,q}\left( x\right)
=G_{n,1}\left( x\right) :=G_{n}\left( x\right) $, where $G_{n}\left(
x\right) $ are known as Genocchi polynomials which is given via the
exponential generating function, as follows:%
\begin{equation*}
\sum_{n=0}^{\infty }G_{n}\left( x\right) \frac{t^{n}}{n!}=\frac{2t}{e^{t}+1}%
e^{xt}\text{ }\left( \left\vert t\right\vert <\pi \right) ,
\end{equation*}

(see [1-4, 12, 13, 21]). From the above, as $q\rightarrow 1$ in (\ref%
{equation 111}), we discover the following corollary:

\begin{corollary}
The following%
\begin{gather*}
\sum_{k=0}^{n}\frac{\binom{n}{k}}{\left( k+2\right) \left( k+1\right) }%
G_{k+2}\left( x\right) \\
=\frac{2x^{n+1}}{n+1}-\frac{G_{n+2}\left( x\right) }{\left( n+1\right)
\left( n+2\right) }-\frac{1}{n+1}
\end{gather*}%
is true.
\end{corollary}

Let us take $y=1$ and $n\rightarrow 2n$ into (\ref{equation 105}), becomes%
\begin{align}
& \sum_{k=0}^{n}\frac{\binom{2n}{2k}}{\left( k+1\right) \left( 2k+1\right) }%
G_{2k+2,q}\left( x\right)  \label{equation 112} \\
& =\frac{G_{2n+2,q}\left( x-1\right) +G_{2n+2,q}\left( x+1\right) }{\left(
2n+1\right) \left( 2n+2\right) }  \notag \\
& =\frac{\frac{1}{q}\left( qG_{2n+2,q}\left( x+1\right) +G_{2n+2,q}\left(
x\right) \right) +qG_{2n+2,q}\left( x\right) +G_{2n+2,q}\left( x-1\right) }{%
\left( 2n+1\right) \left( 2n+2\right) }  \notag \\
& -\frac{G_{2n+2,q}\left( x\right) }{q\left( 2n+1\right) \left( 2n+2\right) }%
-\frac{qG_{2n+2,q}\left( x\right) }{\left( 2n+1\right) \left( 2n+2\right) } 
\notag \\
& =\frac{2\left( n+2\right) x^{n+1}}{\left( 2n+1\right) \left( 2n+2\right) }+%
\frac{2\left( n+2\right) \left( x-1\right) ^{n+1}}{\left( 2n+1\right) \left(
2n+2\right) }-\frac{G_{2n+2,q}\left( x\right) }{q\left( 2n+1\right) \left(
2n+2\right) }-\frac{qG_{2n+2,q}\left( x\right) }{\left( 2n+1\right) \left(
2n+2\right) }  \notag
\end{align}

After these applications, we conclude with the following theorem:

\begin{theorem}
The following identity%
\begin{gather}
\sum_{k=0}^{n}\frac{\binom{2n}{2k}}{\left( k+1\right) \left( 2k+1\right) }%
G_{2k+2,q}\left( x\right)  \label{equation 113} \\
=\frac{\left( n+2\right) x^{n+1}}{\left( 2n+1\right) \left( n+1\right) }+%
\frac{\left( n+2\right) \left( x-1\right) ^{n+1}}{\left( 2n+1\right) \left(
n+1\right) }-\frac{G_{2n+2,q}\left( x\right) }{q\left( 2n+1\right) \left(
2n+2\right) }-\frac{qG_{2n+2,q}\left( x\right) }{\left( 2n+1\right) \left(
2n+2\right) }  \notag
\end{gather}%
is true.
\end{theorem}

Now, we analyse as $q\rightarrow 1$ for the equation (\ref{equation 113})
and so we readily state the following corollary which seems to be
interesting property for the Genocchi polynomials.

\begin{corollary}
The following equality holds:%
\begin{gather}
\sum_{k=0}^{n}\frac{\binom{2n}{2k}}{\left( k+1\right) \left( 2k+1\right) }%
G_{2k+2}\left( x\right)  \label{equation 114} \\
=\frac{2\left( n+2\right) x^{n+1}}{\left( 2n+1\right) \left( 2n+2\right) }+%
\frac{2\left( n+2\right) \left( x-1\right) ^{n+1}}{\left( 2n+1\right) \left(
2n+2\right) }-\frac{G_{2n+2}\left( x\right) }{\left( 2n+1\right) \left(
2n+2\right) }-\frac{G_{2n+2}\left( x\right) }{\left( 2n+1\right) \left(
2n+2\right) }\text{.}  \notag
\end{gather}
\end{corollary}

Substituting $n\rightarrow 2n+1$ and $y=1$ into (\ref{equation 107}), we
compute%
\begin{align*}
& \sum_{k=0}^{n}\frac{\binom{2n+1}{2k}}{\left( k+1\right) \left( 2k+1\right) 
}G_{2k+2,q}\left( x\right) \\
& =\frac{G_{2n+3,q}\left( x+1\right) -G_{2n+3,q}\left( x-1\right) }{\left(
2n+3\right) \left( 2n+2\right) } \\
& =\frac{\frac{1}{q}\left( qG_{2n+3,q}\left( x+1\right) +G_{2n+3,q}\left(
x\right) \right) -\left( G_{2n+3,q}\left( x\right) +G_{2n+3,q}\left(
x-1\right) \right) }{\left( 2n+3\right) \left( 2n+2\right) } \\
& +\left( \frac{q-1}{q}\right) \frac{G_{2n+3,q}\left( x\right) }{\left(
2n+3\right) \left( 2n+2\right) } \\
& =\frac{x^{2n+2}}{q\left( n+1\right) }-\frac{\left( x-1\right) ^{2n+2}}{%
\left( n+1\right) }+\left( \frac{q-1}{q}\right) \frac{G_{2n+3,q}\left(
x\right) }{\left( 2n+3\right) \left( 2n+2\right) }\text{.}
\end{align*}

Therefore, we obtain the following theorem:

\begin{theorem}
The following equality holds:%
\begin{gather*}
\sum_{k=0}^{n}\frac{\binom{2n+1}{2k}}{\left( k+1\right) \left( 2k+1\right) }%
G_{2k+2,q}\left( x\right) \\
=\frac{x^{2n+2}}{q\left( n+1\right) }-\frac{\left( x-1\right) ^{2n+2}}{%
\left( n+1\right) }+\left( \frac{q-1}{q}\right) \frac{G_{2n+3,q}\left(
x\right) }{\left( 2n+3\right) \left( n+1\right) }\text{.}
\end{gather*}
\end{theorem}

As $q\rightarrow 1$ in the above theorem, then we easily derive the
following corollary:

\begin{corollary}
For $n\in 
\mathbb{N}
^{\ast }$, then we have%
\begin{gather*}
\sum_{k=0}^{n}\frac{\binom{2n+1}{2k}}{\left( k+1\right) \left( 2k+1\right) }%
G_{2k+2,q}\left( x\right) \\
=\frac{2x^{2n+2}}{\left( 2n+2\right) }-\frac{2\left( x-1\right) ^{2n+2}}{%
\left( 2n+2\right) }\text{.}
\end{gather*}
\end{corollary}

\section{\textbf{CONCLUSION}}

In this final section, we recall the generating function of the $q$-Genocchi
polynomials, as follows:%
\begin{equation}
\mathcal{F}_{q}\left( x,t\right) =\frac{2t}{qe^{t}+1}e^{xt}=\sum_{n=0}^{%
\infty }G_{n,q}\left( x\right) \frac{t^{n}}{n!}\text{.}  \label{equation 115}
\end{equation}

Using the definition of the $k$-th derivative as $\frac{d^{k}}{dt^{k}}$ to (%
\ref{equation 115}), then we easily see that%
\begin{equation}
\frac{d^{k}}{dt^{k}}\left( \frac{2t}{qe^{t}+1}e^{xt}\right) =\frac{d^{k}}{%
dt^{k}}\left( \sum_{n=0}^{\infty }G_{n,q}\left( x\right) \frac{t^{n}}{n!}%
\right) \text{.}  \label{equation 116}
\end{equation}

By applying $\lim_{t\rightarrow 0}$ on the both sides in (\ref{equation 116}%
), then we conclude with the following theorem:

\begin{theorem}
The following identity%
\begin{equation}
G_{k,q}\left( x\right) =\lim_{t\rightarrow 0}\left[ \frac{d^{k}}{dt^{k}}%
\left( \frac{2t}{qe^{t}+1}e^{xt}\right) \right]  \label{equation 117}
\end{equation}%
is true.
\end{theorem}

We now consider Cauchy-integral formula of the $q$-Genocchi polynomials
which is a vital and important in complex analysis, is an important
statement about line integrals for holomorphic functions in the complex
plane. So, by using the equation of (\ref{equation 117}), we can state the
following theorem:

\begin{theorem}
The following Cauchy-integral holds true:%
\begin{equation*}
G_{n,q}\left( x\right) =\frac{n!}{2\pi i}\int_{C}\mathcal{F}_{q}\left(
x,t\right) \frac{dt}{t^{n+1}}
\end{equation*}%
where $C$ is a loop which starts at $-\infty $, encircles the origin once in
the positive direction, and the returns $-\infty $.
\end{theorem}

Distribution formula for the $q$-Genocchi polynomials is important to study
regarding $p$-adic Measure theory. That is,%
\begin{eqnarray*}
\int_{%
\mathbb{Z}
_{p}}q^{y}\left( x+y\right) ^{n}d\mu _{-1}\left( y\right)
&=&\lim_{n\rightarrow \infty }\sum_{\xi =0}^{dp^{n}-1}\left( -1\right) ^{\xi
}\left( x+\xi \right) ^{n}q^{\xi } \\
&=&d^{n}\sum_{a=0}^{d-1}\left( -1\right) ^{a}q^{a}\left( \lim_{n\rightarrow
\infty }\sum_{\xi =0}^{p^{n}-1}\left( -1\right) ^{\xi }\left( \frac{x+a}{d}%
+\xi \right) ^{n}q^{d\xi }\right) \\
&=&d^{n}\sum_{a=0}^{d-1}\left( -1\right) ^{a}q^{a}\frac{G_{n+1,q}\left( 
\frac{x+a}{d}\right) }{n+1}.
\end{eqnarray*}

After the above applications, we procure the following theorem.

\begin{theorem}
For $n\in 
\mathbb{N}
^{\ast }$, then we have%
\begin{equation*}
G_{n,q}\left( dx\right) =d^{n-1}\sum_{a=0}^{d-1}\left( -1\right)
^{a}q^{a}G_{n,q}\left( x+\frac{a}{d}\right) \text{.}
\end{equation*}
\end{theorem}

By utilizing from the definition of the geometric series in (\ref{equation
115}), we easily see that%
\begin{equation*}
\sum_{m=0}^{\infty }G_{m,q}\left( x\right) \frac{t^{m}}{m!}%
=\sum_{m=0}^{\infty }\left( 2\left( m+1\right) \sum_{n=0}^{\infty }\left(
-1\right) ^{n}q^{n}n^{m}\right) \frac{t^{m+1}}{\left( m+1\right) !},
\end{equation*}

by comparing the coefficients on the both sides, then we have, for $m\in 
\mathbb{N}
$%
\begin{equation}
\frac{G_{m+1,q}\left( x\right) }{m+1}=\left[ 2\right] _{q}\sum_{n=1}^{\infty
}\left( -1\right) ^{n}q^{n}n^{m}\text{.}  \label{equation 118}
\end{equation}

Now also, we develop the following applications by using Mellin
transformation to the generating function of the $q$-Genocchi polynomaials:
For $s\in 
\mathbb{C}
$ and $\Re e\left( s\right) >1$,%
\begin{eqnarray}
\zeta \left( s,x:q\right) &=&\frac{1}{\Gamma \left( s\right) }%
\int_{0}^{\infty }t^{s-2}\left\{ -\mathcal{F}_{q}\left( x,-t\right) \right\}
dt  \label{equation 119} \\
&=&2\sum_{n=0}^{\infty }\left( -1\right) ^{n}q^{n}\left( \int_{0}^{\infty
}t^{s-1}e^{-nt}dt\right)  \notag \\
&=&2\sum_{n=1}^{\infty }\frac{\left( -1\right) ^{n}q^{n}}{n^{s}}  \notag
\end{eqnarray}

Thus, we can state the definition of the $q$-Zeta-type function as follows:%
\begin{equation}
\zeta \left( s,x:q\right) =2\sum_{n=1}^{\infty }\frac{\left( -1\right)
^{n}q^{n}}{n^{s}}\text{.}  \label{equation 120}
\end{equation}%
As $q\rightarrow 1$ in (\ref{equation 120}), turns into%
\begin{equation*}
\lim_{q\rightarrow 1}\zeta \left( s,x:q\right) =\zeta \left( s,x:1\right)
:=\zeta \left( s,x\right) =2\sum_{n=1}^{\infty }\frac{\left( -1\right) ^{n}}{%
n^{s}}
\end{equation*}

which is well-known as Euler-Zeta function (see \cite{KIM}). By (\ref%
{equation 118}) and (\ref{equation 120}), we get 
\begin{equation*}
\zeta \left( -n,x:q\right) =\frac{G_{n+1,q}\left( x\right) }{n+1}\text{.}
\end{equation*}


\begin{thebibliography}{99}
\bibitem{Araci 1} S. Araci, D. Erdal, and J. J. Seo, A study on the
fermionic $p$-adic $q$-integral representation on $%
\mathbb{Z}
_{p}$ associated with weighted $q$-Bernstein and $q$-Genocchi polynomials, 
\textit{Abstract and Applied Analysis}, Volume \textbf{2011}, Article ID
649248, 10 pages.

\bibitem{Araci 2} S. Araci, J. J. Seo, and D. Erdal, New Construction
weighted ($h,q$)-Genocchi numbers and polynomials related to Zeta type
function, \textit{Discrete Dynamics in Nature and Society}, Volume \textbf{%
2011}, Article ID 487490, 7 pages.

\bibitem{araci 3} S. Araci, M. Acikgoz, and Feng Qi, On the $q$-Genocchi
numbers and polynomials with weight zero and their applications, Available
at http://arxiv.org/abs/1202.2643

\bibitem{Araci 4} H. Jolany and H. Sharifi, Some results for the
Apostol-Genocchi polynomials of higher order, \textit{Bulletin of Malaysian
Mathematical Sciences Society} (Article in press).

\bibitem{He} Yuan He and Chunping Wang, Some Formulae of Products of the
Apostol-Bernoulli and Apostol-Euler Polynomials, \textit{Discrete Dynamics
in Nature and Society}, vol. \textbf{2012}, Article ID 927953, 11 pages.

\bibitem{KIM} T. Kim, Euler numbers and polynomials associated with Euler
Zeta functions, \textit{Abstract and Applied Analysis,} Volume \textbf{2008}%
, Article ID 581582, 11 pages.

\bibitem{Kim} T. Kim and S.-H. Lee, Some Properties on the $q$-Euler Numbers
and Polynomials, \textit{Abstract and Applied Analysis}, vol. \textbf{2012},
Article ID 284826, 9 pages, 2012.

\bibitem{Kim 1} D. Kim, T. Kim, S. -H. Lee, D. V. Dolgy, S.-H. Rim, Some new
identities on the Bernoulli numbers and polynomials, \textit{Discrete
Dynamics in Nature and Society}, Volume \textbf{2011,} Article ID 856132, 11
pages.

\bibitem{Kim 2} T. Kim, J. Choi, and Y. H. Kim, Some identities on the $q$%
-Bernoulli numbers and polynomials with weight $0,$ \textit{Abstract And
Applied Analysis}, Volume \textbf{2011}, Article ID 361484, 8 pages.

\bibitem{Kim 3} T. Kim, On a$\ q$-analogue of the $p$-adic log gamma
functions related integrals, \textit{J. Number Theory}, 76 (\textbf{1999})
no. 2, 320-329.

\bibitem{Kim 4} T. Kim, and J. Choi, On the $q$-Euler numbers and
polynomials with weight $0$, \textit{Abstract and Applied Analysis}, Volume 
\textbf{2012}, Article ID 795304, 7 pages..

\bibitem{Kim 5} T. Kim, On the $q$-extension of Euler and Genocchi numbers, 
\textit{J. Math. Anal. Appl.} 326 (\textbf{2007}) 1458-1465.

\bibitem{Kim 6} T. Kim, On the multiple $q$-Genocchi and Euler numbers, 
\textit{Russian J. Math. Phys.} 15 (4) (\textbf{2008}) 481-486.

\bibitem{Kim 7} T. Kim, On the weighted $q$-Bernoulli numbers and
polynomials, \textit{Advanced Studies in Contemporary Mathematics} 21(%
\textbf{2011}), no.2, p. 207-215.

\bibitem{Kim 8} T. Kim, $q$-Volkenborn integration, \textit{Russ. J. Math.
phys. }9 (\textbf{2002}) 288-299.

\bibitem{Kim 9} T. Kim, $q$-Euler numbers and polynomials associated with $p$%
-adic $q$-integrals, \textit{J. Nonlinear Math. Phys.}, 14 (\textbf{2007}),
no. 1, 15--27.

\bibitem{Kim 10} T. Kim, New approach to $q$-Euler polynomials of higher
order, \textit{Russ. J. Math. Phys.,} 17 (\textbf{2010}), no. 2, 218--225.

\bibitem{Kim 11} T. Kim, Some identities on the $q$-Euler polynomials of
higher order and $q$-Stirling numbers by the fermionic $p$-adic integral on $%
\mathbb{Z}
_{p}$, \textit{Russ. J. Math. Phys.,} 16 (\textbf{2009}), no.4,484--491.

\bibitem{Kim 12} D. Kim, T. Kim, J. Choi, and Y. H. Kim, Identities
involving $q$-Bernoulli and $q$-Euler numbers, \textit{Abstract and Applied
Analysis,} Volume \textbf{2012} (2012), Article ID 674210, 10 pages.

\bibitem{C. S. Ryoo} C. S. Ryoo, A note on the weighted $q$-Euler numbers
and polynomials, \textit{Advan. Stud. Contemp. Math. }21(\textbf{2011}),
47-54.

\bibitem{Cangul} I. N. Cangul, V. Kurt, H. Ozden, and Y. Simsek, On the
higher-order $w$-$q$-genocchi numbers, \textit{Advanced Studies in
Contemporary Mathematics }19 (\textbf{2009}), no. 1, pp. 39--57.

\bibitem{Acikgoz} M. Acikgoz, S. Araci and I. N. Cangul, A note on the
modified $q$-Bernstein polynomials for functions of several variables and
their reflections on $q$-Volkenborn integration, \textit{Applied Mathematics
and Computation}, vol. 218 (\textbf{2011}), no. 3, pp. 707--712.

\bibitem{Andrews} L. C. Andrews, Special Functions of Mathematics for
Engineerings, \textit{SPIE Press}, 1992, pages 479.
\end{thebibliography}
\end{document}